\documentclass[a4paper]{article}

\usepackage{framed,multirow}

\usepackage{amssymb}
\usepackage{latexsym}

\usepackage{url}
\usepackage{xcolor}
\definecolor{newcolor}{rgb}{.8,.349,.1}

\usepackage{hyperref}

\usepackage{textcomp}
\usepackage{longtable}
\usepackage{caption}
\usepackage{amssymb}
\usepackage{amsmath}

\usepackage{listings}

\usepackage{paralist}

\usepackage{array}
\usepackage{booktabs}

\usepackage{url}
\usepackage{psfrag}
\usepackage{color}
\usepackage{etoolbox}
\usepackage[detect-all]{siunitx}
\usepackage{booktabs}
\usepackage{etoolbox}
\usepackage{mathtools}
\usepackage{calc}

\robustify\itshape 
\robustify\upshape 
\robustify\mdseries 
\robustify\slshape 

\DeclarePairedDelimiter{\norm}{\lVert}{\rVert}

\newcommand{\ssp}{\hspace*{0.0833333em}} % half of /,
\newcommand{\HN}{{\nu \ssp N}}
\newcommand{\Ml}{{M \ssp l}}
\newcommand{\ark}{a}

\newcommand{\BB}{\beta^{\HN}_M}
\newcommand{\BBone}{\beta^{\HN}_1}

\newcommand{\AL}{\alpha^{\HN}_M}
\newcommand{\ALZ}{\alpha^{\HN}_{a \ssp M}}
\newcommand{\ALE}{\alpha^{\HN}_{s \ssp M}}

\newcommand{\secrow}{\slshape}

\newcommand{\csos}[1]{} % blind comments

\newcommand{\mli}[1]{{\mathit{#1}}}
\newcommand{\safe}{\mli{safe}}

\newcommand{\rmsz}[1]{\mbox{{\sz{\rm{#1}}}}}

\newcommand{\mymat}[1]{{#1}}
\newcommand{\beq}{\begin{equation}}
\newcommand{\eeq}{\end{equation}}
\newcommand{\bec}{\begin{center}}
\newcommand{\eec}{\end{center}}
\newcommand{\pone}[2]{\frac{{\partial} #1}{\partial {#2}}}
\newcommand{\ptwo}[2]{\frac{{\partial}^2 #1}{\partial {#2}^2}}

\newcommand{\n}{\noindent}

\newcommand{\TT}{T}

\newcommand{\ii}{{\rm{i}}}

\newcommand{\sq}{\hspace*{0.5 ex}}

\newcommand{\gesim}{\,\raisebox{-0.4ex}{$\stackrel{>}{\scriptstyle\sim}$}\,}
\newcommand{\lesim}{\,\raisebox{-0.4ex}{$\stackrel{<}{\scriptstyle\sim}$}\,}

\newcommand{\sz}[1]{{\scriptsize#1}}
\newcommand{\RKG}{\sz{RKG} }
\newcommand{\FRKG}{\sz{FRKG} }
\newcommand{\FRKC}{\sz{FRKC} }

\newcommand{\SRK}{\sz{SRK} }
\newcommand{\RK}{\sz{RK} }

\newcommand{\ROCK}{\sz{ROCK} }

\newcommand{\RKLone}{\sz{RKL1} }
\newcommand{\RKLtwo}{\sz{RKL2} }
\newcommand{\RKC}{\sz{RKC} }

\newcommand{\DUMKA}{\sz{DUMKA} }

\newcommand{\pFRKG}{\sz{FRKG}}
\newcommand{\pFRKGcontrol}{\sz{FRKGc}}
\newcommand{\pFRKCcontrol}{\sz{FRKCc}}
\newcommand{\pFRKGfixed}{\sz{FRKGf}}
\newcommand{\pFRKCfixed}{\sz{FRKCf}}

\newcommand{\pRKG}{\sz{RKG}}
\newcommand{\pSRK}{\sz{SRK}}

\renewcommand*{\Re}{\operatorname{Re}} 
\renewcommand*{\Im}{\operatorname{Im}}

\makeatletter % http://www.matthewflickinger.com/blog/archives/2005/02/20/latex_mod_spacing.asp
\def\imod#1{\allowbreak\mkern10mu({\operator@font mod}\,\,#1)}
\makeatother

\newcommand{\reva}[1]{#1}
\newcommand{\revb}[1]{#1}
\newcommand{\revbbox}[1]{#1}

\newcommand{\myv}{\rule{0pt}{3ex}}

\usepackage{cleveref}
\usepackage{amsmath}
\usepackage{amsthm}
\usepackage{algorithm}
\usepackage{algorithmic}
\usepackage{subfigure}
\usepackage{lineno}

\numberwithin{theorem}{section}
\usepackage{amsopn}

\DeclareMathOperator{\erf}{erf}

\usepackage{authblk}
\newcommand{\sep}{, }
\providecommand{\keywords}[1]
{
  \noindent
  \small	
  \textbf{\textit{Keywords:}} #1
}
\providecommand{\msccodes}[1]
               {
                 \noindent
  \small	
  \textbf{\textit{2010 MSC:}} #1
}

\begin{document}

%\begin{frontmatter}

\title{Runge--Kutta--Gegenbauer explicit methods for advection-diffusion problems}%

\author[1]{Stephen O'Sullivan\thanks{Email address: \url{stephen.osullivan@dit.ie}}}

\affil[1]{School of Mathematical Sciences, Technological University Dublin, Kevin Street, Dublin 8, Ireland}

\maketitle

\begin{abstract}
%%%

  In this paper, Runge--Kutta--Gegenbauer (\pRKG) stability polynomials of arbitrarily high order of accuracy are introduced in closed form. The stability domain of \RKG polynomials extends in the the real direction with the square of polynomial degree, and in the imaginary direction as an increasing function of Gegenbauer parameter. Consequently, the polynomials are naturally suited to the construction of high order stabilized Runge--Kutta (\pSRK) explicit methods for systems of PDEs of mixed hyperbolic-parabolic type.
  
  We present \SRK methods composed of $L$ ordered forward Euler stages, with complex-valued stepsizes derived from the roots of \RKG stability polynomials of degree $L$. Internal stability is maintained at large stage number through an ordering algorithm which limits internal amplification factors to $10 L^2$. Test results for mildly stiff nonlinear advection-diffusion-reaction problems with moderate ($\lesim 1$) mesh P\'eclet numbers are provided at second, fourth, and sixth orders, with nonlinear reaction terms treated by complex splitting techniques above second order.

%%%%
\end{abstract}


\msccodes{65L04 \sep 65L20 \sep  65M20}

%\end{frontmatter}

%\linenumbers


\section{Introduction}
\label{intro}

Stabilized Runge--Kutta (\pSRK) explicit methods are particularly well suited to solving mildly stiff systems of ODEs arising from the discretization of parabolic PDEs due to their extended stability domains along the negative real axis. In this work, we will consider an extension of this class of methods to systems of ODEs derived from PDEs of mixed hyperbolic-parabolic type.

\revb{The canonical $m$-dimensional scalar advection-diffusion equation for a quantity $w$ is given by
    \beq
    w_t+\sum_{k=1}^ma_kw_{x_k}=d\sum_{k=1}^mw_{x_k x_k} ,
    \label{eqn:refprob}
    \eeq
    where $a_k$ is the advection coefficient in the $k$-th direction and $d$ is the diffusion coefficient.}
Assuming a spatial mesh of uniform spacing $h_k$ in the $k$-th direction, spatial discretization leads to a system of ODEs via the method of lines which may be written in the form
\beq
w'=f(t,\,w) .%,\quad t>0,\quad w(0)=w_0 
\label{eqn:firsteqn}
\eeq
\revb{Given an initial state $w_0$, subsequent application of a numerical integration method  over $n$ steps then yields an approximate solution $w^n$ at time $t^n$.}

Under the description given above, the mesh P\'eclet number associated with the $k$-th dimension, $P_k=\lvert a_k\rvert h_k/d$, describes the relative significance of the hyperbolic to the parabolic parts of \cref{eqn:refprob}. In particular, the domain of the eigenvalues from a von Neumann stability analysis increases in the imaginary direction with $P_k$. A well constructed numerical method must therefore capture this stability domain to avoid \revb{unbounded growth of errors}.

\SRK methods based on Chebyshev polynomials (eg. \RKC~\cite{sommeijer1997rkc}, \DUMKA~\cite{medovikov1998high}, \ROCK~\cite{abdulle2002fourth}, and \FRKC~\cite{OSULLIVAN2015665}) are suitable for problems in the limit of vanishing P\'eclet numbers. The stability domains associated with these polynomials are extended along the real axis, however, in unmodified form, they also possess internal points where the domains become vanishingly narrow. In practice, damping procedures have been used to introduce a finite imaginary extent at these points of marginal stability so as to mitigate against instability arising through truncation errors. A damping process has been exploited in \cite{verwer2004rkc} to extend applicability of \RKC methods to problems with moderate ($\lesim 1$) mesh P\'eclet numbers.\footnote{The authors also consider large ($> 1$) mesh P\'eclet numbers, however, the benefits of using \SRK methods in terms of efficiency are largely lost.}  %In principle, \RKL methods are also amenable to this treatment via appropriate damping.
The Extrapolated Stabilized Explicit Runge–Kutta (ESERK) methods of~\cite{MARTINVAQUERO2016141} are derived from Richardson extrapolation techniques and demonstrate finite extent stability domains in the imaginary sense through damping.  Finally,  we remark that methods based on Legendre polynomials have also been proposed~\cite{meyer2012second,meyer2014stabilized} which do not suffer from marginally stable internal points and therefore do not require damping for problems with very small ($\ll 1$) mesh P\'eclet numbers.

In this work, we demonstrate that the necessity for modification of \SRK methods through a damping procedure for problems with mild-to-moderate advection may be avoided by appealing to the properties of the general class of Gegenbauer polynomials. We seek a closed-form prescription for arbitrarily high order Runge--Kutta--Gegenbauer (RKG) stability polynomials which natively generate stability domains with imaginary extent determined by the Gegenbauer parameter.

The paper is organized as follows. In \cref{sec:rkg}, the analytic form of the class of \RKG  stability polynomials is presented and the construction of stable time-marching explicit methods based on the roots of these polynomials is outlined. 
In \cref{sec:tests}, numerical tests are presented confirming the order and efficiency properties of \RKG methods. Conclusions are presented in \cref{sec:conc}.

%http://homepages.cwi.nl/~willem/Coll_AdvDiffReac/notes.pdf The message simply is that the eigenvalue criterion (3.7) should be handled with great care if the matrix is not normal. For a thorough discussion we refer to Morton (1980), Dorsselear et al. (1993). Further stability results can also be found in Strikwerda (1989), Thom´ee (1990).

\section{Runge-Kutta-Gegenbauer methods}
\label{sec:rkg}

\subsection{Runge-Kutta-Gegenbauer stability polynomials}
\label{sec:rkgstab}

By appending $t$ to the vector of dependent variables, \cref{eqn:firsteqn} may be written in autonomous form,
\beq
w'=f(w) . %,\quad t>0,\quad w(0)=w_0
\label{eqn:autoneqn} 
\eeq
% Parentheses may be used in the remainder of this work to differentiate exponents from indices.

\revb{We consider advancing the approximate solution explicitly over a single timestep $\TT$ from $\mymat{w}^{n}$, at time level $n$, to $\mymat{w}^{n+1}$, at time level $n+1$. In the following discussion, an order $N$ \SRK method is implemented over $L$ stages spanning an aggregate time $\TT$. For clarity of notation, we omit the time level indexing and denote the $L+1$ internal stage states $W^l$ ($l=0,\,\cdots,\,L$), with  $\mymat{W}^{0}=w^n$, and $\mymat{W}^{L}=w^{n+1}$.  The \SRK integration then takes the form
\beq
\mymat{W}^{L}=\mymat{W}^0+\TT\sum^{L}_{l=1}\ark_l f(W^{l-1}) ,
\label{eqn:rk}
\eeq
\n where the timestep related to each stage is given by $\tau_l= \ark_l\TT$.
}

    \revb{
      When applied to the scalar Dahlquist test equation
      \beq
      \mymat{W}'(t)=\lambda\mymat{W}(t) , \quad \mymat{W}(0)=1 ,
      \label{eqn:dahl}
      \eeq
      the \RKG method of rank $N$, degree $L=MN$, and Gegenbauer parameter $\nu$, may be written as
      \beq
      \mymat{W}^{L}=R^{\HN}_M(z), 
      \eeq
      where $z=\TT\lambda$.
      The associated \RKG stability function, 
      \beq
      \label{eqn:RKGstab}
      %R^{\HN}_M(z)=G^{\HN}_M\left(1+z/\DD\right) ,
      R^{\HN}_M(z)=G^{\HN}_M\left(1+\frac{2z}{\BB}\right) ,
      \eeq
      is obtained from the shifted \RKG polynomial
      \beq
      \label{eqn:RKG}
      G^{\HN}_M(z) = d^{\HN}_0 + 2\sum_{k=1}^{N} d^{\HN}_k C^\nu_{kM}(z) ,
      \eeq
      \n where $C^\nu_{kM}$ denotes the the Gegenbauer polynomial of degree $kM$ and Gegenbauer parameter $\nu$.~\footnote{\revb{We remark that for $\nu=0$ and $\nu=1/2$ the components are Chebyshev and Legendre polynomials respectively.}}
      We note that since any perturbation will be amplified by $R^{\HN}_M(z)$ over a timestep, the method's stability domain is defined by $\{z\in\mathbb{C}:| R^{\HN}_M(z)|\le 1\}$.
    }
    
    \revb{
      Linear order conditions are obtained by requiring that the first $N$ terms of the Taylor series for the exact solution of \cref{eqn:dahl} and $R^{\HN}_M(z)$, coincide. Hence, the order coefficients, $d^{\HN}_k$, are determined by
      \beq
      R^{\HN}_M{}^{(n)}(0)=1 , \quad n=1,\,\ldots,\,N ,
      \label{eqn:linearorderconds}
      \eeq
      where a superscript $(n)$ indicates the $n$th derivative is taken.
    }
\revb{The order conditions are met by solving an $N$-dimensional linear system}
\beq
\label{eqn-pattern}
\left[ \begin{array}{ccc}
C_{M}^{\nu\,(1)}(1) & \ldots & C_{NM}^{\nu\,(1)}(1) \\
\vdots    & \ddots & \vdots \\
C_{M}^{\nu\,(N)}(1) & \ldots & C_{NM}^{\nu\,(N)}(1) 
\end{array}\right]
\left[\begin{array}{c}
d^{\HN}_1\\
\vdots\\
d^{\HN}_N
\end{array}\right]
=
\frac{1}{2}
\left[\begin{array}{c}
(\BB/2)^1\\
\vdots\\
(\BB/2)^N
\end{array}\right] ,
\eeq

\n coupled with a zeroth order condition

\beq
d^{\HN}_0=1-2\sum_{k=1}^{N} d^{\HN}_k C^\nu_{kM}(1) .
\eeq

\n Once the value of $\BB$ has been set, the above equations are sufficient to determine $R^{\HN}_M$:  for odd $M$, optimal values of $\BB$ are found by solving $G^{\HN}_M(-1)=(-1)^N$ iteratively;\csos{ via bisection SEE polyfinder.c: omitting details of initial values for beta (too messy)}
for even $M$, rational interpolation/extrapolation is used.% to assign values for $\BB$.% (the result is reduced by 0.1\% for safety at small $M$ and $\nu$).

In this work, we consider \RKG polynomials of order $N=1,\,\ldots,\,8$ and adopt geometric sequences of values for $\nu$, from 0 to $2N$, given by $\nu=2^{i/2}N/128$, for $i=0,\,\ldots,\,16$, and also include the reference value $\nu=0$. In considering the gain in efficiency achievable by methods constructed from the \RKG polynomial, it is useful to note that $M$ aggregated steps of a standard \RK methods at order $N$ will have a stability polynomial of degree $L$ and a domain with a real extent of approximately $M\beta^{\HN}_1$. From \cref{fig:betavals}, it may be seen that the real extent, $\BB$, of an \RKG polynomial of the same degree is characteristically longer by a factor of between $M$ and $1.5M$ for $\nu=0$, falling to $\lesim 0.5M$ with increasing $\nu$. This factor therefore represents a significant potential gain in efficiency with respect to standard explicit \RK methods. A count of the number of evaluations, $f_{\rm count}$, for the function $f$ in \cref{eqn:rk}, may be used as a measure of the work required to carry out an integration. We note here that the optimal efficiency for \SRK methods follows $f_{\rm count}\propto ({\rm err})^{-1/2N}$ \revb{(where ${\rm err}$ is a measure of the error in the solution)}.

  \begin{figure}[!tbhp]
  \centering
  % From /home/sdo/Code-2015/c/frkglib-5.1/Stabilitydomain-Reduction/
  \revbbox{
    \resizebox{\columnwidth}{!}{\input{betafig.tex}}
    }
  %\resizebox{0.8\columnwidth}{!}{\input{stabdomainfig.tex}}
  \caption{Stability domain extent, for $\lvert R^{\HN}_M\rvert=1$, along the real axis normalized to $M^2\BBone$ as a function of $M$ for orders $N=1,\,\ldots,\,8$. The Gegenbauer parameter $\nu$ ranges from 0 to $2N$ (from top to bottom) with $\nu=0$ and $\nu=2^{i/2}N/128$, for $i=0,\,\ldots,\,16$.\csos{The visible noise component is an effect of the normalization and is not of practical significance.}}
    \label{fig:betavals}
\end{figure}

\subsection{Factorized Runge--Kutta--Gegenbauer methods}

The Factorized Runge--Kutta--Gegenbauer (\pFRKG) coefficients of order $N$ corresponding to \cref{eqn:rk} are defined according to
\beq
\ark^{\HN}_\Ml=\frac{2}{\BB}~\frac{1}{1-\zeta^{\HN}_\Ml} ,
\label{eqn-taunodamp}
\eeq
where the roots $\zeta^{\HN}_\Ml$ of the \RKG polynomial $G^{\HN}_M(z)$ are determined numerically.\footnote{Root-finding is handled with the MPSolve~\cite{bini2000design} package. Multiple-precision calculations of polynomial and method coefficients are carried out using the GMP~\cite{granlund2004gnu} and MPFR~\cite{fousse2007mpfr} libraries.} The maximum stable stepsize for the derived method is $\TT_{\rm max}=\BB/\lvert \lambda\rvert_{\rmsz{max}}$, where the values of $\lambda$ are the negative-definite eigenvalues of the Jacobian associated with \cref{eqn:autoneqn}. 

Over the extent of the stability domain on the real axis, with $\nu>0$, we remark that $\lvert R^{\HN}_M{}^{(n)}\rvert <1$. \reva{In the limiting case of $\nu=0$,} the components of $R^{\HN}_M(z)$ are shifted Chebyshev polynomials and $M-1$ marginally stable internal points exist for which $\lvert R^{\HN}_M{}^{(n)}\rvert =1$. For problems with formally real eigenvalues under von Neumann stability analysis, truncation errors may give rise to eigenvalues with small imaginary parts \reva{thereby introducing some susceptibility to instability.} In this work, we therefore choose a minimum value of $\nu=N/128$ for all test cases so as to avoid this phenomenon.

For various standard discretization schemes, \cite{wesseling2009principles} determined geometric shapes within which the eigenvalues derived from a von Neumann stability analysis will be contained. While \cite{verwer2004rkc} adopt an oval approximation, we find that better results are found for large $M$ using the ellipse approximation 
\beq
\label{eqn:ellipse}
\left(\frac{x}{\BB/2}+1\right)^2+\left(\frac{y}{\AL}\right)^2=1 ,
\eeq
where the centre point is $(-\BB/2,\,0)$, and the major and minor half-axes are $\BB/2$ and $\AL$ respectively.

\begin{figure}[!tbhp]
  \centering
  % From /home/sdo/Code-2015/c/frkglib-5.1/Stabilitydomain-Reduction/
  %  \resizebox{0.8\columnwidth}{!}{\input{stabdomainfig.tex}}
  \resizebox{\textwidth}{!}{\input{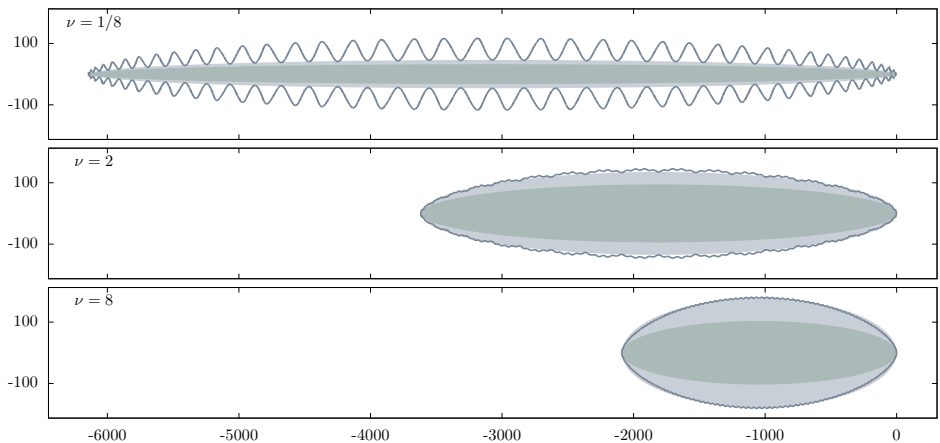}}
  \caption{Stability domains for fourth order \RKG stability polynomials with $M=40$ and $\nu=1/8,\, 2,\, 8$. The shaded regions indicate the strict and approximate fitted ellipses, corresponding to the semi-minor axes $\ALE$ and $\ALZ$ respectively. (The cases shown are characteristic of the methods used in obtaining solutions with $err\gesim 10^{-5}$ for the Brusselator problem with advection presented in \cref{sec:tests} with weak, mild, and moderate advection, from top to bottom.)}
    \label{fig:stabdomain}
\end{figure}

Under $\kappa$-scheme discretizations~\cite{van1985upwind}, the stability condition derived by~\cite{wesseling2009principles} for the ellipse approximation is given by 
\beq
\label{eqn:Tlim}
\TT\le\min\left(
\psi_1\BB ,\,
\psi_2\frac{(\AL)^2}{\BB}
\right) ,
\eeq
where
\beq
\psi_1=\frac{1}{2d\sum h_k^{-2}(2+(1-\kappa)P_k)},\quad\psi_2=\frac{4d}{(2-\kappa)^2\sum a_k^2} .
\label{eqn:psis}
\eeq

%% &M1,&b1,&ae1,&ao1,&az1,&x1,&at0,&at1,&at2,&at3);
%#elif ALPHA==5 /* Midway for fixed, full for stepsize control */
%     a1=(lellipse?ae[M]:ao[M]);
%     at0=at[0][M];
%     redfac=(fixedM>0)?0.5:1.0;  
%     alpha[iGnu][M]=((at0<a1)?a1:(a1+redfac*(at0-a1)));

We identify an ellipse with semi-minor axis, $\ALE$, which is fit strictly to the interior of the contour $\lvert R^{\HN}_M(z) \rvert=1$. Typically, this strict ellipse turns out to be excessively conservative and we therefore also make use of an approximate ellipse with semi-minor axis $\ALZ$. The value of $\ALZ$ is set by solving $\lvert R^{\HN}_M(\ii \ALZ)\rvert =1$ for $M\bmod{4}=0$, and logarithmically interpolating for the intermediate cases.~\footnote{Interpolation yields more conservative values than solving the given equation for $M\bmod{4}\ne 0$.}
%see extraas.c
Except for $M\le 4$, the approximate ellipse semi-minor axis $\ALZ$ is larger than the strict value $\ALE$.

Examples of the stability domains for various \RKG stability polynomials are shown in \cref{fig:stabdomain}, superimposed with the exact and approximate ellipses corresponding to  $\ALE$ and $\ALZ$ respectively. Evidently, the approximate ellipse is well captured within the true stability domain, only overshooting marginally as $\nu$ becomes large. Furthermore, as the amplitude of the ripples in the stability domain boundary decreases with increasing $\nu$, the approximate ellipse becomes an increasingly precise global approximation.

\begin{algorithm}[!tbhp]
\caption{Selection of $M$ and $\nu$ for timestep $\TT$ ($\nu= \nu(i)$, increasing with $i$)}
\label{alg:method}
\begin{algorithmic}
\STATE{Initialize $M=1$,\,$i=1$}
\LOOP
\IF{$\psi_1\BB\le\psi_2 {(\AL)^2}/{\BB}$} % viable width, check extent
\IF{$\BB/\lvert \lambda\rvert _{\rmsz{max}}<\TT$} % viable length: done
\STATE{Select $M$, $\nu= \nu(i)$}  % (with $\nu=\nu(i)$)
\ELSE
\STATE{Increment $M$}
\ENDIF
\ELSE
\STATE{Increment $i$} %\COMMENT{Increasing $\nu$}
\ENDIF
\ENDLOOP
\end{algorithmic}
\end{algorithm}

\begin{figure}[!tbhp]
  \centering
  % From /home/sdo/Code-2015/c/frkglib-5.1/Stabilitydomain-Reduction/
  \revbbox{
    \resizebox{\columnwidth}{!}{\input{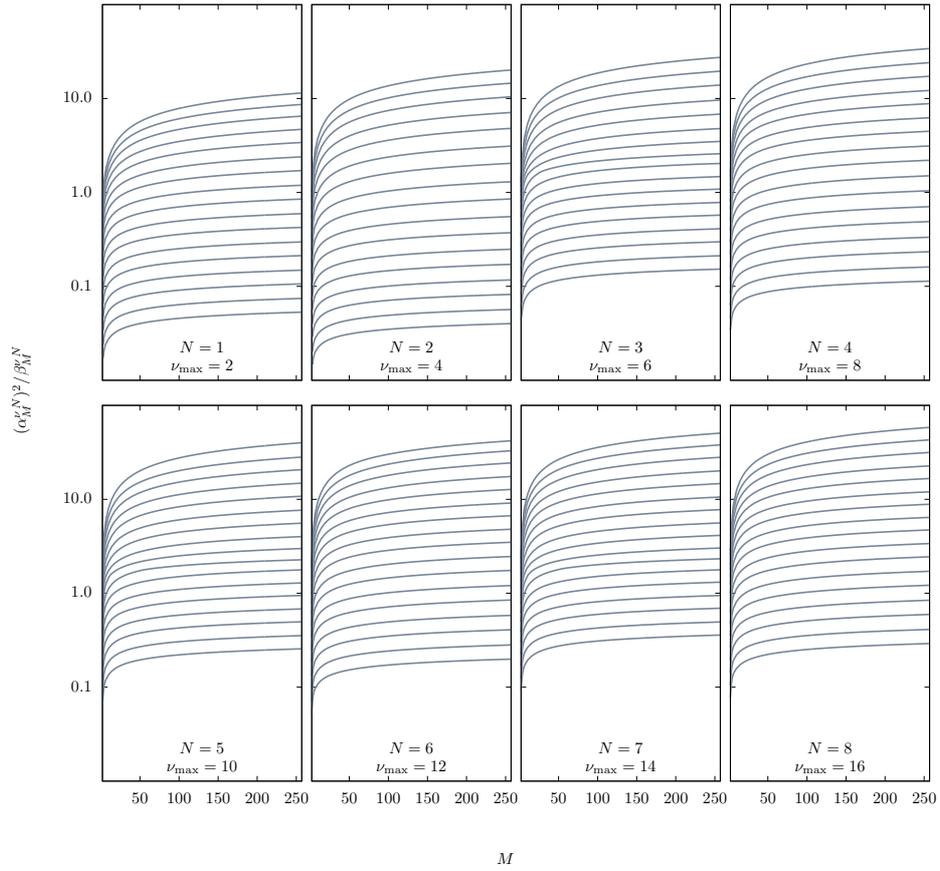}}
  }
  %\resizebox{0.8\columnwidth}{!}{\input{stabdomainfig.tex}}
  \caption{The parameter $(\AL)^2/\BB$ governing the stable timestep limit for the ellipse approximation to the stability domain as a function of $M$ for orders $N=1,\,\ldots,\,8$ with $\nu=2^{i/2}N/128$, for $i=0,\,\ldots,\,16$ (bottom to top). As shown, $(\AL)^2/\BB$ is approximately linear as a function of the Gegenbauer parameter $\nu$.}
  \label{fig:az1vals}
\end{figure}

In cases where a fixed timestep is assigned, we find experimentally that $\AL=\max(\ALE,\,(\ALE+\ALZ)/2)$ is an effective choice for \cref{eqn:Tlim}. Alternatively, when a timestep controller is employed (discussed further in \cref{sec:control}), integration is less sensitive to capturing the stability domain completely and $\AL=\max(\ALE,\,\ALZ)$ proves to be a suitable choice. (The $\max$ function evaluates to $\ALE$ only in a limited number of cases for $M\le 4$.) As illustrated in \cref{fig:az1vals}, the value of $(\AL)^2/\BB$ varies approximately linearly with $\nu$ over two orders of magnitude (from $\sim 0.1$ to $\sim 10$) for $0 \lesim \nu \le 2N$, and also demonstrates a general upward trend with increasing order $N$. The integration method is chosen according to \cref{alg:method} by selecting from the minimum available values of $M$ and $\nu$ for which $\psi_1\BB\le\psi_2(\AL)^2/\BB$.

\subsection{Internal stability}
\label{sec:intstab}

Internal instability arises when the product of any internal sequence of stages generates large values which drown out available numerical precision. Following an idea of Lebedev~\cite{lebedev2000explicit,lebedev1994solve}, but using the more effective method presented in~\cite{OSULLIVAN2015665}, and described here in \cref{alg:intstab}, we order the stages in the $L$-tuple $[\ark^{\HN}_\Ml]$ to approximately minimize the maximum realized internal amplification factor 
 \beq
 \mathcal{Q}=\max(\mathcal{Q}_{j,\,k}(x)), \quad 1\le j,\,k\le L, \quad x\in[-\BB,\,0] , 
 \eeq
 where
 \beq
 \mathcal{Q}_{j,\,k}(x)=\prod_{l=j}^k\lvert 1+\ark^{\HN}_\Ml x\rvert .
 \eeq

 This is a rapid calculation: trialing each swap consists of multiplying the pre-existing forward amplification factor by $v_{j,\,k}$, and dividing the reverse factor by the same quantity. Adopting a slightly reduced stability range, $[(1-10^{-3} n)\BB,\,0]$ ($n\in\mathbb{N}$), aids in meeting the imposed upper bound of $10 L^{2}$ (particularly for small values of $\nu$ and $N$, where the mean value of $n$ reaches almost $3$). Preserving 7 digits for precision, a method consisting of $\sim 10^4 $ stages is therefore theoretically viable in a numerical integration carried out at double (16 digit) precision - well beyond the requirements of practical calculations.

\begin{figure}[!tbhp]
  \centering
  % From /home/sdo/Code-2015/c/frkglib-5.1/Stabilitydomain-Reduction/
  %  \resizebox{0.8\columnwidth}{!}{\input{stabdomainfig.tex}}
    \resizebox{\textwidth}{!}{\input{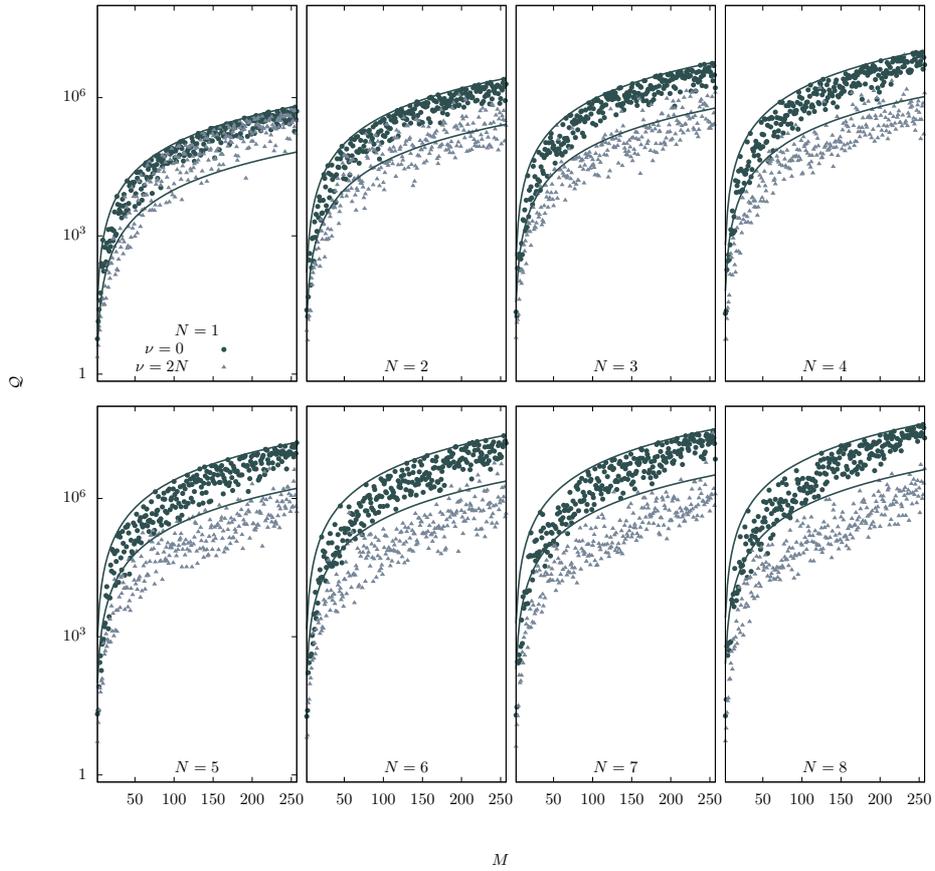}}
\caption{The maximum realized internal amplification factor $\mathcal{Q}$ as a function of $M$ for \FRKG integration methods over orders $N=1,\,\ldots,\,8$. The number of internal stages is given by $L=MN$. Cases for $\nu=0$ (filled circles) and $\nu=2N$ (filled triangles) are shown with guidelines at $L^2$ and $10L^2$.}
\label{fig:intstabfig}
\end{figure}

The maximum realized internal amplification factor $\mathcal{Q}$ is illustrated for $\nu=0$ and $\nu=2N$ over orders $N=1,\,\ldots,\,8$ in \cref{fig:intstabfig}. The imposed bound $10 L^{2}$ is observed in all cases, with $\mathcal{Q}$ falling away from this limit as $\nu$ increases.\footnote{\reva{The ordered coefficients from first to eighth order, with $M=1,\,\ldots,\,257$, are provided as Electronic Supplementary Material.}}

 % SEE ordering13.c!!!
\begin{algorithm}[!tbhp]
  \caption{Internal stability ordering of $L$-tuple $[\ark^{\HN}_\Ml]$}
  \label{alg:intstab}
  \begin{algorithmic}
    \STATE{Shuffle $[\ark^{\HN}_\Ml]$}
    \STATE{Initialize $n=1$}
    \REPEAT
    \STATE{$L$ points $x_k$ uniformly on $[(1-10^{-3}n)\BB,\,0]$, with $v_{j,\,k}=\lvert 1+\ark_j x_k\rvert $}
    \FOR{$l=1$ \TO $l=L$ } % DETERMINE WHICH TERM GOES IN POSITION l
    \FOR{$m=l$ \TO $m=L$ } % Trial \ark^{\HN}_{M\ssp m} in pos l
    \STATE{Swap positions of $\ark^{\HN}_\Ml$ and $\ark^{\HN}_{M\ssp m}$}
    \IF{
      $\norm*{\max\left(\prod_{j=1}^{l}v_{j,\,k},\,\prod_{j=l+1}^{L}v_{j,\,k}\right)}_1$
      is not new minimum over $m$
    }
    \STATE{Revert positions of swapped coefficients $\ark^{\HN}_\Ml$ and $\ark^{\HN}_{M\ssp m}$}
    \ENDIF
    \ENDFOR
    \ENDFOR
    \STATE{Increment $n$}
    \UNTIL{$\mathcal{Q}<10 L^{2}$}
    \COMMENT{confirmed over $10L$ uniformly spaced points}
  \end{algorithmic}
\end{algorithm}

\subsection{Complex splitting}
\label{sec:cmplx}

Up to second order, \FRKG methods are suitable without further consideration for nonlinear problems (since all order conditions are linear). At higher orders, nonlinear order conditions are present which require additional treatment. One approach to ensuring these nonlinear conditions are met is the composition of \RK methods~\cite{zbMATH05292586,hairernonstiff}, which offers a means of combining a linear \SRK method with finishing stages derived from a nonlinear \RK method~\cite{medovikov1998high,abdulle2002fourth}. While composition methods are elegant, our studies and \cite{hundsdorfer2003numerical} suggest that poor error propagation properties over the finishing stages gives rise to erratic behaviour. Furthermore, the latter authors assert that composition methods are prone to order reduction whereby temporal accuracy is lost when spatial mesh refinement occurs.

\n An alternative approach for a linear evolution equation of the form $w'=(A+B)w$, is to employ a high order approximation of the operator ${\rm e}^{\TT(A+B)}$ in terms of complex valued coefficients~\cite{Castella09,hansenostermann09}. In this work, we adopt splitting of the form
\beq
\mymat{w}^{n+1}={\rm e}^{\TT_{k_J}  \mymat{B}}{\rm e}^{\TT_{k_{J-1}}  \mymat{A}}\cdots{\rm e}^{\TT_{k_3}  \mymat{B}}{\rm e}^{\TT_{k_2}  \mymat{A}}{\rm e}^{\TT_{k_1}  \mymat{B}} \mymat{w}^n ,
\label{eqn:split}
\eeq
\n where $\TT_{k}\in\Re^+$ for even $k$, corresponding to the partial timesteps for the operator $A$, and $\TT_{k}\in\Im^+$ otherwise~\cite{Castella09,blanes2013optimized}. 

The procedure extends immediately to the case $w'=Aw+f_B(w)$, with $B$ replaced by a nonlinear operator. The factors ${\rm e}^{\TT_{k}  \mymat{B}}$ are replaced by appropriate approximations to the true flows given by $f_B(w)$ over the the intervals $\TT_{k}$. 

Under Dirichlet or Neumann boundary conditions, splitting methods are known to be prone to an order reduction phenomenon. However, the effect appears to remain confined to the neighbourhood of the boundaries such that the order may be fully recovered on the interior of the domain~\cite{hansenostermann09,lubich1995interior}. 

\subsection{Stepsize control}
\label{sec:control}

For reliable numerical integration, an appropriate control procedure based on local error estimation is required to manage stepsize selection.

    \revb{Starting from $\mymat{w}^{n}$, the approximate solutions $\mymat{w}^{n+1}$ and $\overline{\mymat{w}}^{n+1}$ are obtained at order $N$, and at some lower order of accuracy $\overline{N}$, respectively. A measure of the error over the step is then given by}

\beq
err^{n+1}=\norm*{\frac{\overline{\mymat{w}}^{n+1}-\mymat{w}^{n+1}}{wt}}_2 ,
\label{eqn:err}
\eeq
where
\beq
wt=atol+\max(|\mymat{w}^{n}|,\,|\mymat{w}^{n+1}|)\times rtol ,
\eeq
with $atol$ and $rtol$ being tuning parameters for the absolute and relative errors respectively.

If $err^{n+1}>1$, the step is rejected and $\mymat{w}^{n+1}$ is recalculated with a revised timestep $\TT_{\rm new}$. Otherwise, the solution is accepted and a trial solution $\mymat{w}^{n+2}$ is determined over $\TT_{\rm new}$. In order to prescribe $\TT_{\rm new}$, it is observed that the error behavior for the \FRKG methods follows $err\approx C \TT^{(N+1)/(N-\overline{N})}$. Hence, the constant $C$ may be specified in terms of the error measure, and the revised stepsize may be estimated. For unsplit problems (when the previous step has been accepted) this is done by means of the the predictive controller~\cite{watts1984,gustafsson1994control,hairerstiff} given by
\beq
\TT^{n+1}_{\rm new}= \safe \TT^n\left(\frac{1}{err^{n+1}}\right)^{\frac{N-\overline{N}}{N+1}}
\left(\frac{\TT^{n}}{\TT^{n-1}}\right)
\left(\frac{err^{n}}{err^{n+1}}\right)^{\frac{N-\overline{N}}{N+1}} .
\label{eqn:errcont}
\eeq
In all other cases, the non-predictive value obtained by deleting the third term in parentheses in \cref{eqn:errcont} is used. In practice, we choose $\overline{N}=1$ for second order integrations, and $\overline{N}=2$ at higher orders where splitting is applied. Additionally, we set the safety factor $\safe$ to 0.8 (after the initial step) and allow the revised timestep $\TT_{\rm new}$ to vary by at most a factor of two with respect to the previously trialed value. 

For the initial step, following~\cite{sommeijer1997rkc}, the error is estimated over a trial step $\TT_{\rm trial}=1/\lvert \lambda\rvert_{\rmsz{max}}$ by comparing forward Euler steps using function evaluations calculated at $t=0$ and $t=\TT_{\rm trial}$. The initial step is then assigned using $\TT=\safe\times \TT_{\rm trial}/\sqrt{err}$, with $\safe=0.1$.

\subsection{Convex Monotone Property}
% See http://gwu.geverstine.com/pdenum.pdf for stencil nomenclature

Methods based on Chebyshev polynomials applied to pure diffusion problems with spatially varying diffusion coefficients may give rise to a staircase profile in the solution when discontinuities are present~\cite{meyer2014stabilized}. The effect arises when the coefficients of the equation stencil admit negative values. Methods which have strictly non-negative stencil coefficients for a given problem are described as possessing the Convex Monotone Property (CMP) by \cite{meyer2014stabilized}. These authors also note that that damping is effective in recovering the CMP in Chebyshev polynomial based methods. However, they assert that to maintain the CMP with increasing stage number $L$, damping must be increased, with a resultant reduction in efficiency such that the extent of the stability domain on the real axis goes as $L^{3/2}$.

In this section, we consider the CMP in the more general case of advection-diffusion problems of the form
\beq
w_t+a w_x= d w_{x x} ,
\label{eqn:advdiff1d}
\eeq
in the presence of an initial discontinuity
\beq
w(x,\,0) =
\begin{cases}
w_L  & \text{if $x<0$,} \\
w_R & \text{otherwise}  
\end{cases} .
\label{eqn:advdiff1dinit}
\eeq
The exact solution to \cref{eqn:advdiff1d,eqn:advdiff1dinit} is given by
\beq
w(x,\,t)=\frac{1}{2}\left[(w_R+w_L)+(w_R-w_L) \erf\left( \frac{x-at}{2\sqrt{dt}} \right)\right] .
\eeq
A uniform mesh is assumed with spacing $h=0.1$ over the domain $-20<x<20$, with Dirichlet boundary conditions imposed using the exact solution. The particular parameters chosen here are $a=0.2$, $d=1$, giving a mesh P\'eclet number of $P=0.2$. We consider a first order upwind discretization of the advection term and a second order centred discretization of the diffusion term.

\begin{figure}[!tbhp]
  \centering
  % From /home/sdo/Code-2015/c/frkglib-5.1/sisc/CMPtests/19Aug
  \resizebox{\columnwidth}{!}{\input{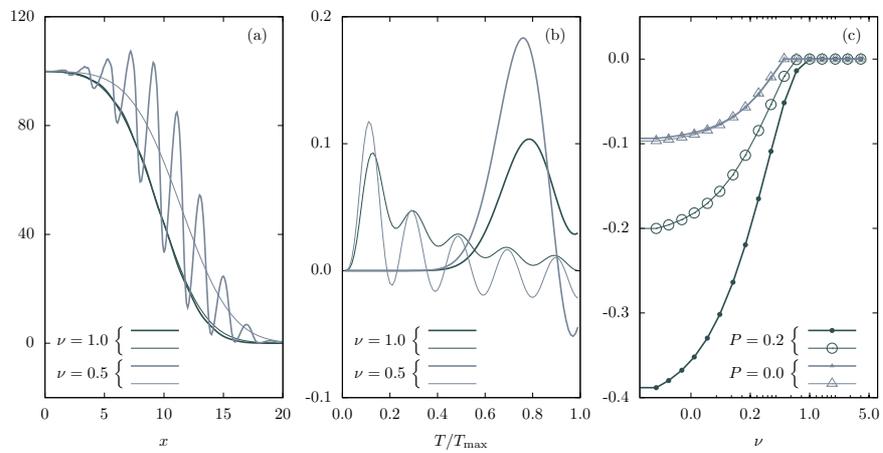}}
  \subfigure{\label{fig:CMPa}} % dummy subfigures for labels only
  \subfigure{\label{fig:CMPb}}
  \subfigure{\label{fig:CMPc}}
  \caption{Illustrations of CMP for advection-diffusion problems. Panel (a): Solutions to \cref{eqn:advdiff1d,eqn:advdiff1dinit} for $\nu=0.5$ and 1.0 at second order accuracy over 10 timesteps for $M=11$ at times 4.9710 and 4.1425 respectively, shown in heavy lines. The exact solutions are also shown in fine lines. Panel (b): Examples of the equation stencil coefficients with timesteps up to the stable limit. For the given case, there are 45 coefficients. Heavy lines show the 39th coefficient (the largest negative value). Fine lines show the 29th coefficient. Panel (c): Mimima of equation stencil coefficients as a function of $\nu$ for $P=0.0$ and $P=0.2$. Heavy lines correspond to $M=11$ and fine lines to $M=7$.
\label{fig:CMP}
  }
\end{figure}

The effect of the CMP is evident in \cref{fig:CMPa}, where solutions from \FRKG methods are presented at second order for $M=11$ over 10 timesteps, with comparisons of the cases $\nu=0.5$ (corresponding to a Legendre polynomial-based method) and $\nu=1$. While both cases are stable, imposing $M> 1$ on the initial discontinuity results in severe oscillations in the solutions derived with $\nu=0.5$, due to an absence of the CMP. We note that the observed oscillatory behaviour is no longer present in the case with $\nu=1$. Evidence as to the origin of this behaviour may be found by inspecting the sample coefficients from the equation stencil shown over the stable range in \cref{fig:CMPb}: at several points, the coefficients corresponding to $\nu=0.5$ become negative, whereas the coefficients for $\nu=1$ are strictly non-negative over the stable range. 

The minimum equation stencil coefficient values are illustrated in \cref{fig:CMPc} for different mesh P\'eclet numbers, $P=0,\,0.2$, and stage numbers $M=7,\,11$. The suggested characteristics of the CMP for \FRKG methods are that it is not maintained for odd $N>1$, however, for even $N$, or $N=1$, the CMP is recovered for some critical value $\nu_{\rm CMP}$. Furthermore, this value of $\nu_{\rm CMP}$ is independent of $M$ for $P=0$, occurring at 0.5 for $N=1,\,2$, and subsequently rising slightly for even $N$. This is consistent with the observed behaviour of the first and second order Legendre polynomial based methods \RKLone and \RKLtwo~\cite{meyer2014stabilized}. In the case of finite advection with $P=0.2$, $\nu_{\rm CMP}$ appears to rise with $M$ more rapidly than the stable estimate provided by \cref{alg:method}. As a consequence, while increasing $\nu$ will significantly dampen the observable staircase artefacts, it is not feasible in practice to fully enforce the CMP for advection-diffusion problems at large $M$ above second order.

We emphasize that the discussions presented here are based on fixed stage number integrations, and that with suitable stepsize control procedures, which will restrict the stage number while discontinuities are diffused, absence of CMP does not generate staircasing in the numerical solution~\cite{1742-6596-837-1-012020}.

\section{Brusselator with advection}
\label{sec:tests}

We consider the Brusselator diffusion-reaction problem~\cite{lefever1971chemical,hairerstiff} extended to include advection,

\begin{eqnarray}
\label{eqn:bruss}
\pone{v}{t}&=&\epsilon\left(\ptwo{v}{x_1} +\ptwo{v}{x_2}\right) +A -(B+1)v +w v^2 +\mu\left(U_1\pone{v}{x_1}+U_2\pone{v}{x_2}\right),\nonumber\\
\pone{w}{t}&=&\epsilon\left(\ptwo{w}{x_1} +\ptwo{w}{x_2}\right)  +B v  -v^2w +\mu\left(V_1\pone{w}{x_1}+V_2\pone{w}{x_2}\right) ,
\end{eqnarray}

%\begin{eqnarray}
%\label{eqn:bruss}
%\plone{v}{t}&=&\epsilon\left(\pltwo{v}{x_1} +\pltwo{v}{x_2}\right) +A -(B+1)v +w v^2 +\mu\left(U_1\plone{v}{x_1}+U_2\plone{v}{x_2}\right),\nonumber\\
%\plone{w}{t}&=&\epsilon\left(\pltwo{w}{x_1} +\pltwo{w}{x_2}\right)  +B v  -v^2w +\mu\left(V_1\plone{w}{x_1}+V_2\plone{w}{x_2}\right) ,
%\end{eqnarray}

\n with initial conditions $v(0,\,x)=22x_2(1-x_2)^{1.5}$, $w(0,\,x)=27x_1(1-x_1)^{1.5}$. Tests are configured with diffusion and reaction parameters $\epsilon=0.01$, $A=1.3$, $B=1$; advection vectors $U=(-0.5,1)^T$, $V=(0.4,0.7)^T$; and $\mu=0.1,\,0.5,\,1.0$, corresponding to weak, mild, and moderate advection respectively. Numerical solutions are obtained at $t=1$ on a uniform mesh with periodic boundary conditions over $0\le x_1\le 1$, $0\le x_2\le 1$, with $h_1=h_2=1/800$. We employ second order upwind discretizations of the advection terms (corresponding to $\kappa=-1$ under the $\kappa$-scheme formalism) with second order centred discretization of the diffusion terms. For the stated advection parameters, the largest mesh P\'eclet numbers (given by $P_2$ for species $v$) are 0.0125, 0.0625, 0.125. The order of the system of ODEs is $2/h^2=640000$, with the spectral radius of the Jacobian matrix estimated by $\psi_1^{-1}\approx 5.17\times 10^4,\,5.36\times 10^4,\,5.60\times 10^4$ respectively. Hence, the problem is moderately stiff and not readily tractable via standard \RK techniques.

\begin{figure}[!tbhp]
  \centering
  % From /home/sdo/Code-2015/c/frkglib-5.1/Stabilitydomain-Reduction/
  \resizebox{\columnwidth}{!}{\input{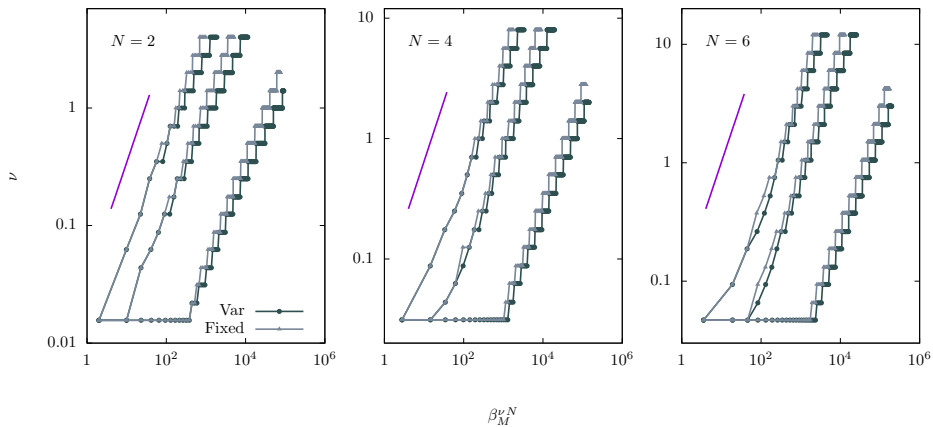}}
  %\resizebox{0.8\columnwidth}{!}{\input{stabdomainfig.tex}}
  \caption{The Gegenbauer parameter, $\nu$, as selected by \cref{alg:method}, plotted against $\BB$ for fixed timestep methods (grey lines) and variable timestep methods using a stepsize controller (black lines). $\nu$ ranges from $N/128$ to $2N$ at 17 logarithmically uniformly spaced values. The methods corresponding to the weak, mild and moderate advection ($P_2=0.0125,\, 0.0625,\, 0.125$) are shown from left to right. Guidelines are linear in $\BB$.}
  \label{fig:method}
\end{figure}

Complex splitting of the form given by \cref{eqn:split} is employed at fourth and sixth orders using coefficients from~\cite{Castella09} and \cite{blanes2013optimized} (also presented in \cref{tab:splitting}). Reaction terms are treated using standard adaptive \RK integrators from the GNU Scientific Library~\cite{galassi2009gnu}. 

\begin{figure}[!tbhp]
  \centering
  % From /home/sdo/Code-2015/c/frkglib-5.1/frkgbrusselator-Reduction/
    \resizebox{\textwidth}{!}{\input{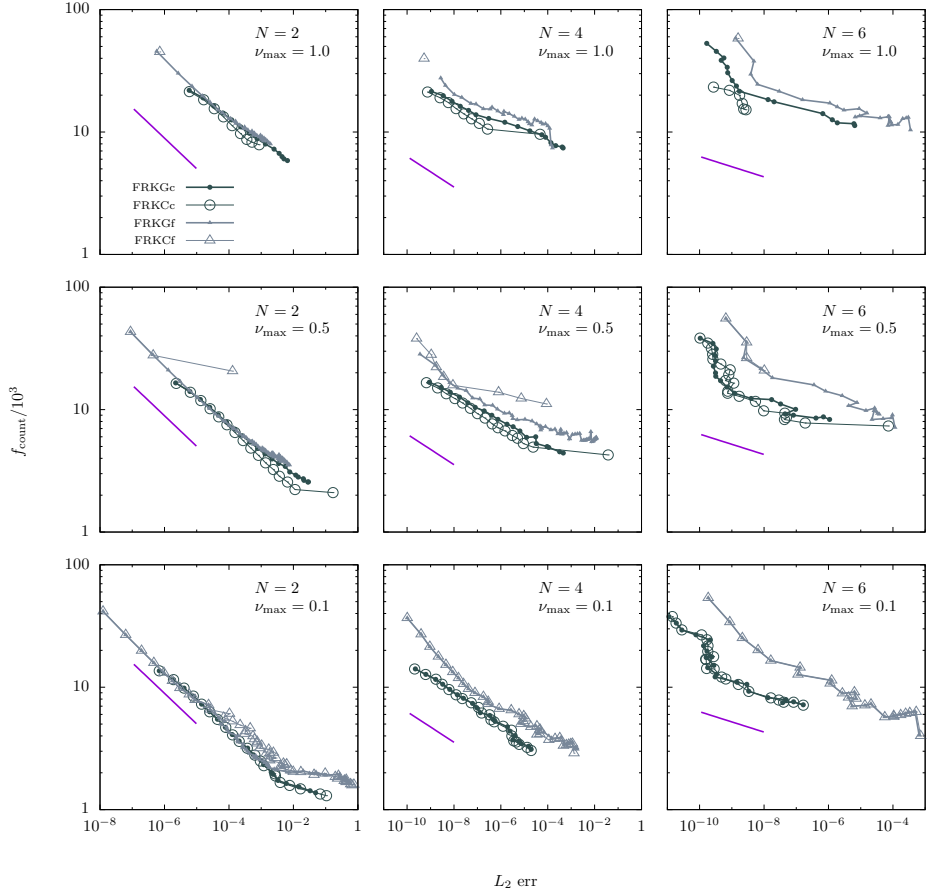}}
    \caption{
      Work-precision diagram for Brusselator problem with advection. Second order unsplit method results are shown in the left column, with fourth and sixth order split method results in the second and third columns. Mesh P\'eclet numbers, from bottom row to top, are 0.0125, 0.0625, 0.125, corresponding to weak, mild, and moderate advection respectively. Data from \FRKG method integrations with stepsize control (\pFRKGcontrol) and fixed stepsizes (\pFRKGfixed) are shown as filled circles and filled triangles respectively. Unfilled circles and triangles correspond to results from Chebyshev polynomial based methods with stepsize control (\pFRKCcontrol) and fixed stepsizes (\pFRKCfixed). Guidelines shown are proportional to $({\rm err})^{-1/2N}$.}
\label{fig:frkgbrusselator}
\end{figure}

The methods implemented in the tests are selected via \cref{alg:method}, for both fixed stepsizes, and stepsizes controlled via local error estimates. The associated Gegenbauer parameter, $\nu$, is shown as a function of the maximum timestep size in \cref{fig:method}. It can be seen that increasing $\nu$ up to the maximum value of $2N$ extends the stable stepsize by approximately two orders of magnitude, whilst continuing to contain the domain of the methods' eigenvalues within the ellipse approximations given by \cref{eqn:ellipse}. The realized timesteps represent a gain in efficiency over standard \RK integrations of $1.5M$, down to $0.5M$ with increasing $\nu$. Typical \RKG stability domains for the integration of the linear parts of \cref{eqn:bruss} encountered at fourth order may be seen in \cref{fig:stabdomain}.

Work-precision test results for the Brusselator problem with advection are shown in \cref{fig:frkgbrusselator} for unsplit second order and split higher order \FRKG methods. Cases are considered with both fixed stepsize over the full integration, and variable stepsizes controlled by the error estimation procedure described in \cref{sec:control}. (Indicative numerical data for the tests are also provided in \cref{app:worktables}.) For comparison, results from integrations carried out with Chebyshev polynomial based \FRKC methods are also presented. 

Stepsize control is clearly an effective counter-measure to the growth of numerical instabilities. In fact, where the advective term is weak-to-mild, or where the internal stage number is low, the stepsize controller is adequate in managing instability growth arising due to eigenvalues with imaginary parts extending into the exterior of the method's stability domain. For large internal stage number with large advection, as seen in the top row of \cref{fig:frkgbrusselator}, instability develops excessively within a single timestep and stepsize control is no longer effective. We note that splitting appears to generate a complex error behaviour which is significantly overestimated by \cref{eqn:err}. This is managed by adopting tolerance parameters $atol$ and $rtol$ for the controller at values larger than the target precision in \cref{eqn:err}. We find that integrations utilizing stepsize control are more efficient in general than fixed stepsize integrations, particularly for higher order split methods.

In the absence of stepsize control, the Chebyshev polynomial based \FRKC methods with weak damping are largely ineffective, with the exception of cases where the advection term is very small (see the bottom row of \cref{fig:frkgbrusselator}). As noted by \cite{verwer2004rkc}, increasing damping improves the performance of these methods. In contrast, \RKG methods rely solely on the natural characteristics of the Gegenbauer polynomials forming the underlying \RKG stability polynomials and do not require explicit damping procedures. 

\section{Conclusions}
\label{sec:conc}

In this paper, we have presented the class of Runge--Kutta--Gegenbauer (\pRKG) stability polynomials in closed form to arbitrarily high order of accuracy. The \RKG polynomials of order $N$, and degree $L=MN$, comprise a linear combination of Gegenbauer polynomials of degree $kM$, for $k=1,\,\ldots,\,N$, and common Gegenbauer parameter. The particular weighting of the combination is chosen to conform to the linear order conditions, subject to maximizing the extent of the stability domain along the negative real axis, which scales as $L^2$. Crucially, for the consideration of systems of mixed hyperbolic-parabolic type, the domain extends in the imaginary direction as an increasing function of the Gegenbauer parameter.

We have demonstrated the construction of Factorized Runge--Kutta--Gegenbauer (\pFRKG) explicit methods to high order, consisting of ordered sequences of forward Euler stages with complex-valued stepsizes. The algorithm implemented in ordering the $L$ stages of a given method prevents internal amplification factors from overwhelming available numerical precision by bounding their magnitudes to $10 L^2$.

\RKG stability polynomials are shown to be effective in the construction of high order explicit methods for mildly stiff advection-diffusion problems with moderate ($\lesim 1$) mesh P\'eclet numbers. % for which standard \RK methods are inefficient 
~\footnote{\reva{An implementation of second order \FRKG methods is available as Electronic Supplementary Material.}}

\section*{Acknowledgments}
The author wishes to acknowledge the DJEI/DES/SFI/HEA Irish Centre for High-End Computing (ICHEC) for the provision of computational facilities and support.

%\newpage
\reva{
  \section*{Data Statement}
Electronic Supplementary Material for this work is licensed under the Creative Commons Attribution-NonCommercial 4.0 International License.
}

%\newpage
\appendix
\section{Splitting methods}
\label{app:splitting}
%
%PARAM TABLE
\begin{longtable}{SSS[table-format=2.42,table-text-alignment=left]}
\caption{Complex operator splitting parameters, corresponding to \cref{eqn:split}, at orders $N=4,\,6$~\cite{Castella09,blanes2013optimized,blanesN6coeffs}. For each quoted value of $N$, the rows show: $j$, the index for a distinct sweep, with $\Re(T_j)$ and $\Im(T_j)$, the real and imaginary components of the associated timescale; $J$ (in the last row), the number of distinct sweeps required; and the sequence of labels $j$ identifying $k_1\,\cdots\,k_J$.
}
\label{tab:splitting}\\
%\begin{tabular}{SSS[table-format=2.42]}
\hline\noalign{\smallskip}\\
\mbox{$N$}   & \mbox{$j$} & \mbox{$\Re(T_j)$}   \\  % fdeltanew from output?
  &  &  \secrow \mbox{$\Im(T_j)$}   \\  % fdeltanew from output?
%\noalign{\smallskip}\cline{2-3}\noalign{\smallskip}
\myv
& \mbox{$J$} &\mbox{$k_1\,\cdots\,k_J$}   \\ 
%1 &  1 & 1 & 10 & 0.05 & 0.45125 \\
% & & & & & \secrow 0.0\\
%  & & 2 & 10 & 0.05 & 0.45125 \\
% & & & & & \secrow 0.0\\
%2  & 1 & 1.0 \\
%  & & \secrow 0.0\\
% &   2 & 0.5 \\
% &  & \secrow 0.0\\
%\noalign{\smallskip}\cline{2-3}\noalign{\smallskip}
% &  3 & \multicolumn{1}{l}{\mbox{\hspace*{1ex} 2\sq 1\sq 2}}\\
\noalign{\smallskip}\hline\noalign{\smallskip}
4   & 1 &      {$1/4$}   \\
 &  & \secrow  {0}     \\
 &   2 &       {1/10}  \\
 &  & \secrow  {\hspace{-1.5ex}$-$1/30} \\
 &   3 &       {4/15}  \\ 
 &  & \secrow  {2/15}  \\
 &   4 &       {4/15}  \\
 &  & \secrow  {\hspace{-1.5ex}$-$1/5}  \\
%\noalign{\smallskip}\cline{2-3}\noalign{\smallskip}
\myv
&   9 & \multicolumn{1}{l}{\mbox{\hspace*{1ex} 2\sq 1\sq 3\sq 1\sq 4\sq 1\sq 3\sq 1\sq 2}}\\
\noalign{\smallskip}\hline\noalign{\smallskip}
6   & 1 & 0.0625 \\
 &  & \secrow 0.0\\
 &   2 & 0.02469487608701806464091086499684224783860 \\
 &  & \secrow  -0.00787479556290687705817157794952694216320 \\
 &   3 & 0.06381347402130269977936630418820014696320 \\
 &  & \secrow 0.03536576103414332780462940464971474181270 \\
 &   4 & 0.06842509403031644197039700782174468405850 \\
 &  & \secrow -0.06226224445074867699533254064444759604610 \\
 &   5 & 0.08804770109226783762699719586940866757720 \\
 &  & \secrow 0.04547387150229870438376254918797742644469 \\
 &   6 & 0.02368961112984706069614191247000936432533 \\
 &  & \secrow 0.00962432606408962405769803529063730666395 \\
 &   7 & 0.04272972238677338220296430057707421855388 \\
 &  & \secrow -0.03399440392395761055408394845784435826499 \\
 &   8 & 0.12233468631684577296042851700196256307880 \\
 &  & \secrow -0.01043585907975251066938082710059054955178 \\
 &   9 & 0.04189843282969388604353685060726223976426 \\
 &  & \secrow 0.06936249263169638427515817430714426213030 \\
 &   10 & 0.04873280421186970815851409293499173568080 \\
 &  & \secrow -0.09051829642972473048855853856612858205130 \\
%\noalign{\smallskip}\cline{2-3}\noalign{\smallskip}
\myv
&  33 & \multicolumn{1}{l}{\mbox{\hspace*{1ex} 2\sq 1\sq 3\sq 1\sq 4\sq 1\sq 5\sq 1\sq 6\sq 1\sq 7\sq 1\sq 8\sq 1\sq 9\sq 1\sq 10\sq 1\sq 9\sq 1\sq 8\sq 1\sq 7\sq 1\sq 6\sq 1\sq 5\sq 1\sq 4\sq 1\sq 3\sq 1\sq 2}}\\
\noalign{\smallskip}\hline
%\end{tabular}
\end{longtable}

\newpage
\section{Work-precision data tables for Brusselator with advection.}
\label{app:worktables}
%\vspace*{-7mm}
\begin{table}[!tbhp]
  \caption{Second order ($N=2$) \FRKG method results for the Brusselator problem with advection detailed in \cref{sec:tests}. Columns show: $L_2$ (and $L_{\infty}$) errors; count of function evaluations, $f_{\rm count}$; number of timesteps attempted, Timesteps; and the maximum value of $M$ required, $M_{\rm max}$. Weak, mild, and moderate advection data ($\mu=0.1,\,0.5,\,1.0$) are presented, from bottom to top respectively. For rows where values of $atol$ are given, the stepsize controller is implemented and $f_{\rm count}$ values refer to accepted (rejected) function evaluations.}
  \label{tab:brussN2}
  \centering
  %\resizebox{10cm}{!}{
  %\input{frkgbrusselatortable}
  \begin{tabular}{lllll} \hline
 $atol$ & $L_2$ ($L_{\infty}$) error &  $f_{\rm count}$ (rej) & Timesteps & $M_{\rm max}$ \\ \hline
\multicolumn{5}{c}{$N=2$}\\\hline
\multicolumn{5}{c}{$\mu=1.0$}\\\hline
$10^{-2}$ & 5.09e-03 (3.06e-02) &  6030 (0) & 52 & 51 \\
$10^{-3}$ & 1.35e-03 (8.61e-03) &  7975 (0) & 154 & 51 \\
$10^{-4}$ & 1.22e-04 (7.71e-04) &  12559 (0) & 533 & 20 \\
$10^{-5}$ & 5.86e-06 (3.70e-05) &  21898 (0) & 1766 & 9 \\
\myv
  & 1.12e-04 (4.84e-04) &  12870 & 322 & 10 \\
  & 1.08e-03 (4.61e-03) &  9120 & 76 & 30 \\
\multicolumn{5}{c}{\myv$\mu=0.5$}\\\hline
$1$ & 2.81e-02 (2.06e-01) &  2561 (0) & 14 & 126 \\
$10^{-1}$ & 1.40e-02 (9.53e-02) &  2814 (0) & 20 & 126 \\
$10^{-2}$ & 3.49e-03 (2.28e-02) &  3628 (0) & 44 & 122 \\
$10^{-3}$ & 4.52e-04 (2.92e-03) &  5423 (0) & 124 & 48 \\
$10^{-4}$ & 4.66e-05 (3.00e-04) &  9065 (0) & 382 & 21 \\
$10^{-5}$ & 2.22e-06 (1.43e-05) &  16525 (0) & 1209 & 11 \\
\myv
  & 3.75e-05 (1.74e-04) &  9230 & 231 & 10 \\
  & 7.99e-04 (3.67e-03) &  5010 & 42 & 30 \\
  & 1.51e-03 (6.85e-03) &  4480 & 28 & 40 \\
  & 2.11e-03 (9.37e-03) &  4400 & 22 & 50 \\
  & 4.29e-03 (1.89e-02) &  3780 & 16 & 60 \\
  & 3.94e-03 (1.74e-02) &  4130 & 15 & 70 \\
  & 6.59e-03 (2.90e-02) &  3520 & 11 & 80\\
\multicolumn{5}{c}{\myv$\mu=0.1$}\\\hline
$1$ & 4.87e-02 (4.00e-01) &  1373 (0) & 14 & 105 \\
$10^{-1}$ & 3.32e-03 (1.55e-02) &  1733 (0) & 21 & 80 \\
$10^{-2}$ & 9.91e-04 (5.46e-03) &  2517 (0) & 42 & 65 \\
$10^{-3}$ & 1.26e-04 (6.48e-04) &  4106 (0) & 105 & 38 \\
$10^{-4}$ & 1.42e-05 (7.12e-05) &  7277 (0) & 307 & 23 \\
$10^{-5}$ & 6.90e-07 (3.44e-06) &  13609 (0) & 950 & 13 \\
\myv
  & 8.91e-06 (4.74e-05) &  7870 & 197 & 10 \\
  & 2.02e-04 (8.44e-04) &  4000 & 50 & 20 \\
  & 6.78e-04 (3.04e-03) &  2760 & 23 & 30 \\
  & 1.78e-03 (8.67e-03) &  2200 & 14 & 40 \\
  & 2.75e-01 (5.07e+00) &  1950 & 10 & 50 \\
\hline
\end{tabular}

\vspace*{-7mm}
  %}
\end{table}

%\newpage
\vspace*{-7mm}
\begin{table}[!tbhp]
  \caption{Fourth and sixth order ($N=4,\,6$) split \FRKG method sample results.}
\label{tab:brussHO}
\centering
%\resizebox{10cm}{!}{
%\input{frkgbrusselatortable}
\begin{tabular}{lllll} \hline
 $atol$ & $L_2$ ($L_{\infty}$) error &  $f_{\rm count}$ (rej) & Timesteps & $M_{\rm max}$ \\ \hline
\multicolumn{5}{c}{$N=4$}\\\hline
\multicolumn{5}{c}{$\mu=1.0$}\\\hline
$10^{-2}$ & 1.37e-04 (6.68e-04) &  8178 (0) & 22 & 51 \\
$10^{-3}$ & 5.24e-06 (3.16e-05) &  10672 (410) & 49 & 49 \\
$10^{-4}$ & 4.18e-08 (2.74e-07) &  14886 (130) & 107 & 46 \\
$10^{-5}$ & 1.09e-09 (8.14e-09) &  21444 (0) & 235 & 24 \\
\myv
  & 8.74e-08 (3.48e-07) &  16900 & 65 & 10 \\
\multicolumn{5}{c}{\myv$\mu=0.5$}\\\hline
$10^{-1}$ & 1.10e-04 (3.37e-04) &  4912 (0) & 16 & 123 \\
$10^{-2}$ & 1.26e-05 (4.49e-05) &  5936 (0) & 24 & 113 \\
$10^{-3}$ & 7.47e-07 (2.75e-06) &  7996 (326) & 45 & 73 \\
$10^{-4}$ & 3.92e-08 (1.44e-07) &  11244 (116) & 85 & 43 \\
$10^{-5}$ & 8.93e-10 (3.27e-09) &  16718 (0) & 177 & 27 \\
\myv
  & 6.21e-08 (2.52e-07) &  12420 & 48 & 10 \\
  & 1.07e-05 (4.55e-05) &  8320 & 15 & 20 \\
  & 7.99e-05 (3.44e-04) &  7260 & 9 & 30 \\
  & 4.93e-04 (1.83e-03) &  6400 & 6 & 40 \\
  & 2.77e-03 (1.11e-02) &  5600 & 4 & 50 \\
  & 6.29e-03 (2.22e-02) &  5760 & 4 & 60 \\
\multicolumn{5}{c}{\myv$\mu=0.1$}\\\hline
$1$ & 1.95e-05 (6.33e-05) &  3108 (0) & 12 & 101 \\
$10^{-1}$ & 3.22e-06 (1.04e-05) &  3714 (0) & 16 & 77 \\
$10^{-2}$ & 5.69e-07 (2.00e-06) &  4842 (426) & 26 & 66 \\
$10^{-3}$ & 9.29e-08 (3.46e-07) &  6500 (310) & 41 & 53 \\
$10^{-4}$ & 6.59e-09 (2.47e-08) &  9486 (114) & 74 & 40 \\
$10^{-5}$ & 2.26e-10 (8.40e-10) &  14024 (82) & 149 & 28 \\
\myv
  & 2.55e-08 (1.08e-07) &  10640 & 41 & 10 \\
  & 1.12e-05 (4.07e-05) &  5600 & 11 & 20 \\
  & 1.45e-04 (6.48e-04) &  3840 & 5 & 30 \\
  & 1.60e-03 (5.10e-03) &  3200 & 3 & 40 \\
\multicolumn{5}{c}{\myv$N=6$}\\\hline
\multicolumn{5}{c}{$\mu=1.0$}\\\hline
$10^{-2}$ & 1.31e-06 (9.22e-06) &  12584 (0) & 18 & 49 \\
$10^{-3}$ & 1.72e-09 (6.92e-09) &  20164 (1334) & 47 & 49 \\
$10^{-4}$ & 7.18e-10 (3.20e-09) &  33140 (644) & 101 & 48 \\
$10^{-5}$ & 1.67e-10 (7.54e-10) &  53002 (0) & 221 & 25 \\
\myv
  & 3.03e-06 (1.62e-05) &  15080 & 9 & 10 \\
  & 3.52e-04 (1.23e-03) &  10380 & 2 & 30 \\
\multicolumn{5}{c}{\myv$\mu=0.5$}\\\hline
$1$ & 1.08e-06 (5.30e-06) &  8308 (0) & 12 & 123 \\
$10^{-1}$ & 4.47e-08 (1.66e-07) &  9224 (0) & 15 & 123 \\
$10^{-2}$ & 3.98e-09 (1.49e-08) &  12328 (0) & 23 & 123 \\
$10^{-3}$ & 3.17e-10 (1.06e-09) &  17534 (1188) & 44 & 79 \\
$10^{-4}$ & 2.95e-10 (1.25e-09) &  26278 (1484) & 85 & 44 \\
\myv
  & 7.68e-06 (3.09e-05) &  11400 & 7 & 10 \\
  & 9.64e-05 (3.53e-04) &  9040 & 3 & 20 \\
\multicolumn{5}{c}{\myv$\mu=0.1$}\\\hline
$1$ & 1.65e-07 (1.00e-06) &  7188 (0) & 12 & 101 \\
$10^{-1}$ & 1.32e-08 (8.79e-08) &  8222 (0) & 15 & 91 \\
$10^{-2}$ & 3.03e-10 (1.27e-09) &  11090 (1062) & 25 & 67 \\
$10^{-3}$ & 1.60e-10 (5.95e-10) &  15694 (1798) & 42 & 59 \\
$10^{-4}$ & 2.14e-10 (8.17e-10) &  22804 (1468) & 74 & 40 \\
$10^{-5}$ & 1.13e-11 (3.81e-11) &  36440 (1250) & 146 & 29 \\
\myv
  & 1.12e-06 (3.72e-06) &  10720 & 7 & 10 \\
  & 1.05e-04 (4.79e-04) &  5880 & 2 & 20 \\
\hline
\end{tabular}

%}
\end{table}
\newpage

\section{Supplementary material}
Supplementary material related to this article can be found online at \newline\url{https://doi.org/10.1016/j.jcp.2019.03.001}.

\bibliographystyle{acm}
\bibliography{poly25}

\end{document}